\documentclass[12pt]{amsart}
\usepackage{amscd, amssymb, latexsym}

\newcommand{\bC}{\mathbb{C}}
\newcommand{\bCP}{\mathbb{CP}}
\newcommand{\bE}{\mathbb{E}}
\newcommand{\bL}{\mathbb{L}}
\newcommand{\bP}{\mathbb{P}}
\newcommand{\bQ}{\mathbb{Q}}
\newcommand{\bZ}{\mathbb{Z}}

\newcommand{\cE}{\mathcal{E}}
\newcommand{\cM}{\mathcal{M}}
\newcommand{\cO}{\mathcal{O}}

\newcommand{\ch}{\mathrm{ch}}
\newcommand{\td}{\mathrm{td}}
\newcommand{\pt}{\mathrm{pt}}
\newcommand{\rank}{\mathrm{rank}}
\newcommand{\Aut}{\mathrm{Aut}}
\newcommand{\Ext}{\mathrm{Ext}}
\newcommand{\Br}{\mathrm{Br}}
\newcommand{\Sym}{\mathrm{Sym}}
\newcommand{\vir}{\mathrm{vir}}
\newcommand{\Def}{\mathrm{Def}}
\newcommand{\Obs}{\mathrm{Obs}}

\newcommand{\vx}{\mathbf{x}}
\newcommand{\vy}{\mathbf{y}}

\newcommand{\la}{\lambda}

\newcommand{\bu}{\bullet}
\newcommand{\Mbar}{\overline{\cM}}

\newtheorem{dummy}{dummy}

\newtheorem{theo}[dummy]{Theorem}

\theoremstyle{definition}

\newtheorem{exam}[dummy]{Example}

\begin{document}

\title{Lectures on the ELSV Formula}
\author{Chiu-Chu Melissa Liu}
\address{Department of Mathematics,
Columbia University, New York, NY 10027, USA}
\email{ccliu@math.columbia.edu}


\begin{abstract}
The ELSV formula, first proved by Ekedahl, Lando, Shapiro, and Vainshtein, relates Hurwitz
numbers to Hodge integrals. Graber and Vakil gave another proof of the ELSV formula
by virtual localization on moduli spaces of stable maps to $\bP^1$, and also explained
how to simplify their proof using moduli spaces of relative stable maps to the pair $(\bP^1,\infty)$.
In this expository article, we explain what the  ELSV formula is and how to prove 
it by virtual localization on moduli spaces of relative stable maps, following Graber-Vakil. 
This note is based on lectures given by the author at Summer School on  
``Geometry of Teichm\"{u}ller Spaces and Moduli Spaces of Curves'' at Center of Mathematical
Sciences, Zhejiang University,  July 14--20, 2008.
\end{abstract}

\maketitle

\tableofcontents

\section{Introduction}

Deligne and Mumford introduced the notation of stable curves and constructed the moduli
space $\Mbar_g$ of genus $g$ stable curves. \cite{DM}. 
The moduli space $\Mbar_{g,n}$ of $n$-pointed genus $g$ stable curves was constructed
by Knudsen-Mumford and Knudsen \cite{KM, Kn2, Kn3}.
Since Mumford's seminal paper \cite{Mu} in the early 1980s, 
the intersection theory of moduli
spaces of stable curves has been studied extensively.
Evaluations of Hodge integrals
\begin{equation}\label{eqn:hodge}
\int_{\Mbar_{g,n}}\psi_1^{j_1}\cdots \psi_n^{j_n} 
\lambda_1^{k_1}\cdots \lambda_g^{k_g}
\end{equation}
are important and difficult problems in this subject. Here $\psi_i$ is
the first Chern class of the line bundle $\bL_i\longrightarrow\Mbar_{g,n}$ 
whose fiber at the moduli point $[(C,x_1,\ldots,x_n)]$ is the cotangent line
$T^*_{x_i}C$ at the $i$-th marked point; $\lambda_i$ is the $i$-th Chern
class of the Hodge bundle $\bE\longrightarrow \Mbar_{g,n}$, which 
is a rank $g$ complex vector bundle whose fiber at $[(C,x_1,\ldots,x_n)]$ is $H^0(C,\omega_C)$,
the space of sections of the dualizing sheaf $\omega_C$ of $C$ (see Section \ref{sec:hodge}
for a review on Hodge integrals).

Using Mumford's Grothendieck-Riemann-Roch calculations in \cite{Mu},  
Faber proved, in \cite{Fa},  that general Hodge integrals can
be uniquely reconstructed from the $\psi$ integrals 
(also known as {\em descendant integrals}):
\begin{equation}\label{eqn:descendant}
\int_{\Mbar_{g,h}}\psi_1^{j_1}\cdots \psi_h^{j_h}. 
\end{equation}
The descendant integrals can be computed recursively by Witten's conjecture
which asserts that the $\psi$ integrals \eqref{eqn:descendant}
satisfy a system of differential equations known as the KdV equations \cite{Wi}.
The KdV equations and the string equation determine all the $\psi$ integrals
\eqref{eqn:descendant} from the initial value $\int_{\Mbar_{0,3}} 1=1$.

The Witten's conjecture was first proved by Kontsevich in \cite{Ko1}. This is
one of the most striking and fundamental result in the intersection
theory of  moduli spaces of stable curves. By now, Witten's conjecture has been
reproved many times (Okounkov-Pandharipande \cite{OP1},
Mirzakhani \cite{Mi},  Kim-Liu \cite{KiL}, Kazarian-Lando \cite{KaL},  
Chen-Li-Liu \cite{CLL}, Kazarian \cite{Ka}, Mulase-Zhang \cite{MuZ} ...).  
The ELSV formula, which relates Hurwitz numbers to Hodge integrals, plays a central
role in several of the above proofs (Okounkov-Pandharipande, Kazarian-Lando,
Kazarian, Mulase-Zhang ...).
The ELSV formula is named after Ekedahl, Lando, Shapiro, and  Vainshtein,
who first proved this formula in \cite{ELSV1, ELSV2}. Later,
Graber and Vakil gave another proof by virtual localization
on moduli spaces $\Mbar_{g,0}(\bP^1,d)$ of genus $g$ degree $d$ stable
maps to $\bP^1$ \cite{GV1}. Fantechi and Pandharipande proved a special
case of the ELSV formula by virtual localization \cite[Theorem 2]{FanP}. 
Okounkov and Pandharipande's paper \cite{OP1} contains a detailed
exposition of the proof of the ELSV formula by virtual
localization on moduli spaces of stable maps to $\bP^1$, following
\cite{FanP} and \cite{GV1}.

In \cite[Section 5]{GV1}, Graber and Vakil explained how their proof 
could be much simplified using moduli spaces of relative stable maps to 
the pair $(\bP^1,\infty)$. When Graber and Vakil wrote their paper \cite{GV1},
moduli spaces of relative stable maps had already been constructed in the symplectic 
category, by Li-Ruan \cite{LR} and by Ionel-Parker \cite{IP1, IP2}.
However, Graber and Vakil needed such moduli spaces in
the algebraic category, with desired properties (proper Deligne-Mumford
stack with perfect obstruction theory, so that the virtual fundamental
class exists, and the virtual localization is applicable). 
Jun Li constructed moduli spaces of relative stable maps
with desired properties in the algebraic category \cite{Li1, Li2}.
In this expository article, we explain what the ELSV formula is  and
how to prove it by virtual localization on moduli spaces
of relative stable maps, following Graber-Vakil \cite{GV1}.

Virtual localization on moduli spaces of relative stable maps can
be used to prove other Hodge integral identities. 
In \cite{LLZ1} (resp. \cite{LLZ2}), K. Liu, J. Zhou and the author used virtual localization on
moduli spaces of relative stable maps to $\bP^1$ relative to $\infty$ (resp.
to the toric blowup of $\bP^2$ at two torus fixed points relative
to the two exceptional divisors) to prove 
the Mari\~{n}o-Vafa formula \cite{MV} (resp. a formula of two-partition Hodge integrals \cite{Zh2}),
which relates certain generating function of Hodge integrals to the colored HOMFLY invariants
of the unknot (resp. the Hopf link).  
(Okounkov and Pandharipande gave another proof of the 
Mari\~no-Vafa formula using virtual localization on moduli spaces of
stable maps to $\bP^1$ \cite{OP2}.) 
The ELSV formula can be obtained by taking certain limit of the Mari\~{n}o-Vafa formula.
See \cite{Liu} for a survey of proofs and applications of the Mari\~{n}o-Vafa formula
and the formula of two-partition Hodge integrals.

We now give an overview of the remainder of this
paper. In Section \ref{sec:hurwitz-hodge}, we recall the definitions
of Hurwitz numbers and Hodge integrals,  and state the ELSV formula.
In Section \ref{sec:equivariant-localization}, we give a brief review of 
equivariant cohomology and localization. In Section \ref{sec:proof},
we interpret Hurwitz numbers as certain relative Gromov-Witten
of the pair $(\bP^1,\infty)$, and derive the ELSV formula
by virtual localization on moduli spaces of relative
stable maps to the pair $(\bP^1,\infty)$.

\subsection*{Acknowledgments}
This note is based on lectures given by the author at 
Summer School on ``Geometry of Teichm\"{u}ller Spaces and 
Moduli Spaces of Curves'' at Center of Mathematical
Sciences, Zhejiang University, July 14--20, 2008. I wish
to thank the organizers Lizhen Ji, Kefeng Liu, and Shing-Tung
Yau for inviting me to give these lectures and to contribute to 
the proceedings.  I wish to thank Jun Li, Kefeng Liu and 
Jian Zhou for the collaboration, which shapes my
current understanding of the materials presented in this note. 
Finally, special thanks go to Lizhen Ji for his kindness
and patience during the preparation of this note.

\section{Hurwitz Numbers and Hodge Integrals}\label{sec:hurwitz-hodge}

In this section, we review the definitions of Hurwitz numbers
and Hodge integrals, and give the precise statement
of the ELSV formula.

\subsection{Hurwitz numbers}\label{sec:hurwitz}
In this subsection, we give a brief review of the geometric
and combinatorial definitions of Hurwitz numbers  which
count ramified covers of $\bP^1$ with a given
ramification type over $\infty\in \bP^1$. 

Let $d$ be a positive integer. A {\em partition} of
$d$ is a sequence of positive integers
$\mu=(\mu_1\geq \mu_2\geq \cdots\geq \mu_h>0)$
such that $\mu_1+\cdots + \mu_h=d$.
The sum of all components of $\mu$, $d$, is called the {\em size} 
of the partition $\mu$, denoted $|\mu|$; 
the number of components in $\mu$, $h$, is called the {\em length} 
of the partition $\mu$, denoted $\ell(\mu)$.

Given a nonnegative integer $g$ and a partition $\mu=(\mu_1\geq \mu_2\geq \cdots \geq \mu_h>0)$
of a positive integer $d$, we consider ramified covers $f:C\to \bP^1$ satisfying the
following conditions:
\begin{enumerate}
\item[(i)] $C$ is a connected compact Riemann surface of genus $g$.
\item[(ii)] $\deg f = d$.
\item[(iii)] $f^{-1}(\infty) =\sum_{i=1}^h \mu_i x_i$ as
Cartier divisors, where $x_1,\ldots,x_h$ are distinct points in $C$.
(So $\infty$ is a critical value of $f$ if $\ell(\mu)<d$.)

\item[(iv)] All other branch points of $f$ (i.e. critical values of $f$) are 
simple. Namely, if $b\in \bP^1-\{\infty\}$ is a critical value of $f$ then
there is a unique critical point $x\in f^{-1}(b)$, and $x$ is
a nondegenerate critical point. So $f^{-1}(b)$ consists of
exactly $d-1$ distinct points. 
\end{enumerate}

Let $b_1,\ldots,b_r$ be the branch points of $f$ in
$\bC=\bP^1-\{\infty\}$. Then the number $r$ is determined
by the genus $g$ and the ramification type $\mu$. To find 
$r$,  let $B=\{ b_1,\ldots, b_r,\infty\}$, and let $C'= C-f^{-1}(B)$. Then
$f|_{C'}: C'\to \bP^1-B$ is an honest covering map of degree $d$. We have
$$
\chi(C')= d\cdot \chi(\bP^1-B)
$$
where $\chi(C') = 2-2g -(d-1)r -h$, $\chi(\bP^1-B) = 1-r$.
We conclude that 
\begin{equation}
r=2g-2+d+h = 2g-2 +|\mu|+\ell(\mu).
\end{equation}

If we fix $r$ distinct points $b_1,\ldots, b_r \in \bP^1-\{\infty\}$ then
there are only finitely many ramified covers $f:C\to \bP^1$ satisfying
(i)--(iv).  Indeed, the domain Riemann surface $C$ is determined 
by the monodromy $\sigma_i$ around $b_i$ which are transpositions
in the permutation group $S_d$ of $\{1,2,\ldots,d\}$. They satisfy
\begin{equation}
\sigma_1  \cdots \sigma_r =\sigma_\infty 
\end{equation}
where $\sigma_\infty$ is the monodromy around $\infty$. Let
$C_\mu\subset S_d$ be the conjugacy class which consists of
products of $h$ disjoint cycles of lengths $\mu_1,\ldots, \mu_h$.
Then $\sigma_\infty \in C_\mu$.

The connected Hurwitz number $H_{g,\mu}$ counts connected ramified covers
$f:C\to \bP^1$ satisfying (i)--(iv), weighted by $(\#\Aut(f))^{-1}$, where
$\Aut(f)$ is the group of automorphisms of the map $f$ and
$\#\Aut(f)$ denotes the cardinality of  the set $\Aut(f)$.
If we fix a $\sigma_\infty\in C_\mu$ then 
\begin{equation}\label{eqn:connectedH}
\begin{aligned}
H_{g,\mu} =&\frac{1}{z_\mu}
\#\Bigl\{(\sigma_1,\ldots, \sigma_r)\mid \sigma_i\textup{ transpositions in }S_d,\
\sigma_1\cdots \sigma_r =\sigma_\infty,\\
& \quad \quad \quad \langle \sigma_1,\ldots,\sigma_r\rangle
\textup{ acts transitively on} \{1,2,\ldots,d \} \Bigr \},
\end{aligned}
\end{equation}
where $d=|\mu|$, $r=2g-2+|\mu|+\ell(\mu)$, and
$$
z_\mu=\frac{d!}{\# C_\mu}=\mu_1 \cdots\mu_h\cdot \# \Aut(\mu)
$$
is the cardinality of the centralizer of $\sigma_\infty$.

The disconnected Hurwitz number $H^\bullet_{\chi,\mu}$ counts
possibly disconnected ramified covers $f:C\to \bP^1$ satisfying
(i)' $\chi(C)=\chi$ and (ii), (iii), (iv), weighted by $(\#\Aut(f))^{-1}$.
If we fix a $\sigma_\infty\in C_\mu$ then
\begin{equation}\label{eqn:disconnectedH}
H_{\chi,\mu}^\bu =\frac{1}{z_\mu}
\# \Bigl\{(\sigma_1,\ldots, \sigma_r)\mid \sigma_i\textup{ transpositions in }S_d,\
\sigma_1\cdots \sigma_r =\sigma_\infty \Bigr\},
\end{equation}
where $d=|\mu|$ and $r=-\chi+|\mu|+\ell(\mu)$.

Given a partition $\mu$, 
We introduce a generating function $\Phi_\mu(\la)$ (resp.
$\Phi_\mu^\bu(\la)$) of connected (resp. disconnected) Hurwitz numbers
$$
\Phi_\mu(\la) = \sum_{g=0}^\infty \la^{2g-2+\ell(\mu)} H_{g,\mu},\quad
\Phi_\mu^\bu(\la)=  \sum_{\chi} \la^{-\chi+\ell(\mu)} H_{\chi,\mu}^\bu.
$$
We now introduce variables $x_1, x_2,\ldots$ and
let $p_i = x_1^i+ x_2^i+\cdots $ be the Newton polynomials.
Given a partition
$\mu=(\mu_1\geq  \cdots\geq \mu_h>0)$, define
$p_\mu = p_{\mu_1}\cdots p_{\mu_h}$. Then
$$
\exp\Bigl(\sum_{\mu\neq \emptyset} \Phi_\mu(\la)p_\mu\Bigr)
= \sum_\mu\Phi^\bu_\mu(\la)p_\mu.
$$
where $\emptyset$ denotes the empty partition
(the unique partition with zero size and zero length).

Identities \eqref{eqn:connectedH} and \eqref{eqn:disconnectedH}
define Hurwitz numbers in terms of representations
of the permutation group. Given a partition $\mu$ of $d>0$,
let $R_\mu$ denote the irreducible representation of $S_d$ associated to $\mu$, and let 
$\chi_\mu$ be the character of $R_\mu$.  
\begin{theo}[Burnside formula] \label{thm:burnside}
Let $\mu=(\mu_1\geq \mu_2\geq \cdots \geq \mu_h>0)$
be a partition of $d$. Then
$$
\Phi_\mu^\bu(\la)=\sum_{|\nu|=d} \frac{\chi_\nu(C_\mu)}{z_\mu}
e^{\frac{\kappa_\nu \la}{2}} \frac{\dim R_\nu}{d!}.
$$
where $\kappa_\mu = \sum_{i=1}^h \mu_i(\mu_i-2i+1)$.
\end{theo}

Indeed, Theorem \ref{thm:burnside} is a special case of
the Burnside formula for general Hurwitz numbers 
$H^h_{g,\mu^1,\ldots,\mu^k}$, where $g$, $h$ are nonnegative integers,
and $\mu^1,\ldots,\mu^k$  are partitions of the same positive integer $d$.
The Hurwitz number $H^h_{g,\mu^1,\ldots,\mu^k}$ counts, with weight, 
degree $d$ ramified covers $f:C\to D$ of
a fixed genus $h$ Riemann surface $D$ by a genus $g$ Riemann surface $C$,
with prescribed ramification types $\mu^1,\ldots, \mu^k$ over $k$ fixed distinct points
$q_1,\ldots, q_k$ in $D$. Theorem \ref{thm:burnside} corresponds to 
the special case where $h=0$, $k=1$.

\subsection{Hodge Integrals}\label{sec:hodge}
Let $\Mbar_{g,h}$ be the moduli space of $h$-pointed, genus $g$
stable curves. (In this paper we always work over $\bC$.) A point in $\Mbar_{g,h}$ is represented by
$[(C,x_1,\ldots,x_h)]$, where $C$ is a complex algebraic curve of arithmetic genus $g$ with at most
nodal singularities, $x_1,\ldots, x_h$ are distinct smooth points on $C$, and $[(C,x_1,\ldots,x_h)]$ is
{\em stable} in the sense that its automorphism group is finite. When $C$ is smooth, it can
be viewed as a connected compact Riemann surface of genus $g$.

$\Mbar_{g,h}$ is a proper smooth Deligne-Mumford stack (or a compact, complex, smooth orbifold)
of (complex) dimension $3g-3+h$. It has a fundamental class
$[\Mbar_{g,h}]\in H_{2(3g-3+h)}(\Mbar_{g,h};\bQ)$. Given 
$\alpha_1,\ldots,\alpha_k \in H^*(\Mbar_{g,h};\bQ)$, we define
their top intersection number to be
\begin{equation}
\int_{\Mbar_{g,h}}\alpha_1 \cdots  \alpha_k :=
\langle [\Mbar_{g,h}],\alpha_1\cdots \alpha_k\rangle\in \bQ 
\end{equation}
where $\langle \ , \ \rangle$ is the pairing
between the homology $H_*(\Mbar_{g,h};\bQ)$
and the cohomology $H^*(\Mbar_{g,h};\bQ)$.

The dualizing sheaf $\omega_C$ of a curve with at most nodal singularities
is an invertible sheaf (line bundle). Near a smooth point, $\omega_C$ is generated
by the local holomorphic differential $dz$, where $z$ is  a local holomorphic coordinate;
near a node, which is locally isomorphic to $(0,0)\in \{ xy=0\mid (x,y)\in \bC^2\}$,
$\omega_C$ is generated by the meromorphic differential
$dx/x = -dy/y$. The Hodge bundle $\bE$ is a rank $g$ vector bundle over
$\Mbar_{g,h}$ whose fiber over the moduli point
$$
[(C,x_1,\ldots,x_h)]\in \Mbar_{g,h}
$$
is $H^0(C,\omega_C)$.
When $C$ is smooth, $\omega_C =\Omega^1_C$ is the sheaf of
local holomorphic differentials on the compact Riemann surface $C$, 
and $H^0(C,\omega_C)= H^0(C,\Omega^1_C)$
is the space of holomorphic differentials on $C$. 
The $\la$ classes are defined by
$$
\la_j=c_j(\bE)\in H^{2i}(\Mbar_{g,h};\bQ).
$$
The cotangent line $T_{x_i}^* C$ of $C$ at the $i$-th marked point $x_i$
gives rise to a line bundle $\bL_i$ over $\Mbar_{g,h}$.
The $\psi$ classes are defined by
$$
\psi_i=c_1(\bL_i)\in H^2(\Mbar_{g,h};\bQ).
$$
The $\la$ classes and $\psi$ classes lie in $H^*(\Mbar_{g,h};\bQ)$
instead of $H^*(\Mbar_{g,h};\bZ)$ because
$\bE$ and $\bL_i$ are orbibundles on the compact orbifold $\Mbar_{g,h}$.

{\em Hodge integrals} are top intersection numbers of $\lambda$ classes and
$\psi$ classes:
\begin{equation}\label{eqn:hodge}
\int_{\Mbar_{g,h}}\psi_1^{j_1}\cdots \psi_h^{j_h} 
\la_1^{k_1}  \cdots \la_g^{k_g} \in \bQ.
\end{equation}
By definition, \eqref{eqn:hodge} is zero unless
$$
j_1+\cdots + j_r + k_1 + 2k_2 + \cdots + g k_g = 3g-3+h.
$$
A special class of Hodge integrals are {\em linear} Hodge integrals 
\begin{equation}\label{eqn:linear}
\int_{\Mbar_{g,h}}\psi_1^{j_1}\cdots\psi_h^{j_h}\lambda_i
\end{equation}
where $i=0,\ldots,g$. When $i=0$, we have $\lambda_0=1$, 
so \eqref{eqn:linear} reduces to top intersection
of $\psi$ classes, known as {\em descendent integrals}:
$$
\int_{\Mbar_{g,h}}\psi_1^{j_1}\cdots\psi_h^{j_h}.
$$

\subsection{The ELSV formula}
The ELSV formula, first proved by Ekedahl, Lando, Shapiro, and Vainshtein 
\cite{ELSV1, ELSV2}, relates the Hurwitz numbers $H_{g,\mu}$ to linear Hodge integrals.
\begin{theo}[ELSV formula] 
Let $\mu=(\mu_1\geq \mu_2\geq \cdots\geq  \mu_h>0)$ be a partition of $d$. Then
\begin{equation}\label{eqn:ELSV}
H_{g,\mu}=\frac{(2g-2+d+h)! }{\#\Aut(\mu)}
\prod_{i=1}^h\frac{\mu_i^{\mu_i}}{\mu_i!}
\int_{\Mbar_{g,h}}\frac{\Lambda^\vee(1)}{\prod_{i=1}^h(1-\mu_i\psi_i)},
\end{equation}
where 
$$
\Lambda^\vee(1) =\sum_{i=0}^g (-1)^i \lambda_i,\quad
\frac{1}{1-\mu_i\psi_i}= \sum_{j=0}^{3g-3+h}(\mu_i \psi_i)^j.
$$
\end{theo}

\section{Equivariant Cohomology and Localization}
\label{sec:equivariant-localization}

In this section, we give a brief introduction to  equivariant cohomology and
localization. See \cite{Fu} for an excellent (and much more comprehensive)
exposition of this subject.

\subsection{Universal bundle}

Let $G$ be a Lie group. 
Let $EG$ be a contractible topological
space on which $G$ acts freely.
(In this note, all the group actions
are continuous.) 
Suppose that $G$ acts on the
{\em right} of $EG$. 
The quotient $BG= EG/G$ is the classifying 
space of principal $G$-bundles, and  the natural projection
$EG \to BG$ is the universal principal $G$-bundle;
$EG$ and $BG$ are defined up to homotopy equivalences.

\begin{exam}\label{BCstar}
Let $G=\bC^*$. Let $EG = \bC^\infty -\{ 0\}$, which is 
contractible. Let $G=\bC^*$ acts on the right of $EG=\bC^\infty -\{0\}$ by
$$
v\cdot \lambda  = \lambda  v,\quad \lambda\in \bC^*,\quad v\in \bC^\infty-\{0\}.
$$
Then the $G$-action on $EG$ is free. The classifying space
$$
BG=(\bC^n-\{0\})/\bC^*= \bCP^\infty 
$$
is the infinite dimensional complex projective space. 
\end{exam}

In Example \ref{BCstar}, $G=\bC^*$ is a complex algebraic group,  
$EG=\bC^\infty-\{0\}$ is an infinite dimensional complex manifold,
and the $G$-action on $EG$ is holomorphic. So the classifying
space $BG =\bCP^\infty$ is a complex manifold, and
the universal principal $\bC^*$-bundle
$EG\to BG$ is a holomorphic principal $\bC^*$-bundle.
The {\em tautological line bundle} $S\to \bCP^\infty$ is the holomorphic line
bundle associated to the universal principal $\bC^*$-bundle. 
For any $k\in \bZ$, let $\cO_{\bCP^\infty}(k)= S^{\otimes {-k}}$. (Strictly speaking, 
$\cO_{\bCP^\infty}(k)$ is the sheaf of local holomorphic sections
of the holomorphic line bundle $S^{\otimes {-k}}$, but we
will not distinguish $\cO_{\bCP^\infty}(k)$ from $S^{\otimes {-k}}$
in this note.)

If $G=G_1\times G_2$ then we may take $EG= EG_1 \times EG_2$, so 
that $BG = BG_1 \times BG_2$. 

\begin{exam}
If $G=(\bC^*)^n$ then $BG= (B\bC^*)^n = (\bCP^\infty)^n$. 
\end{exam}

\subsection{Equivariant cohomology}
Let $G$ be a Lie group, and let
$X$ be a topological space with a {\em left} $G$-action.
Then $G$ acts on $EG\times X$ freely by 
$$
g\cdot (p,x) = (p\cdot  g^{-1}, g \cdot x).
$$
The {\em homotopy orbit space} $X_G$ is defined
to be the quotient of $EG\times X$ by this free
$G$-action. The projection $EG\times X\to EG$
to the first factor descends to a projection
$\pi: X_G \to BG$, which is a fibration
over $BG$ with fiber $X$.  

The $G$-equivariant cohomology
of the $G$-space $X$ is defined to be the
ordinary cohomology of the homotopy orbit space
$X_G$:
$$
H^*_G(X; R) := H^*(X_G; R)
$$
where $R$ is any coefficient ring.  From now on we
will assume $R=\bQ$, the field of rational numbers,
and write $H^*(\bullet)$ for $H^*(\bullet ;\bQ)$.
The following are some special cases.
\begin{enumerate}
\item If $X$ is point then $X_G =  BG$, so $H_G^*(\pt) = H^*(BG)$.
\item If $G$ acts on $X$ freely then $X_G$ is homotopically equivalent to 
the orbit space $X/G$,  so $H_G^*(X)=H^*(X/G)$.
\item If $G$ acts on $X$ trivially then $X_G = BG \times X$. By K\"{u}nneth formula,
$$
H^*_G(X)\cong H^*(X)\otimes_\bQ H^*(BG).
$$
\end{enumerate}

\begin{exam} If $G=\bC^*$ then $H^*_G(\pt) = H^*(\bCP^\infty) \cong \bQ[u]$,
where  $u\in H^2(X;\bQ)$ is the first Chern class
of $\cO_{\bCP^\infty}(1)$. 
\end{exam}

\begin{exam} If $G=(\bC^*)^n$ then 
$H^*_G(\pt) = H^*((\bCP^\infty)^n)\cong \bQ[u_1,\ldots,u_n]$.
\end{exam}

\begin{exam}\label{Pr}
Let $\bC^*$ act on the $r$-dimensional 
complex projective space $\bP^r$ by 
$$
t\cdot[z_0,\ldots, z_r]= [t^{a_0} z_0,\ldots, t^{a_n} z_n],\quad
t\in \bC^*,\quad [z_0,\ldots, z_r]\in \bP^r,
$$
where $a_0,\ldots, a_r\in \bZ$. 
Then the fibration $\bP^r_{\bC^*} \to B\bC^*$ can be identified
with the $\bP^r$-bundle 
$$
\bP(\cO_{\bCP^\infty}(a_0)\oplus \cdots \oplus \cO_{\bCP^\infty}(a_r)) \to \bCP^\infty.
$$
To compute $H^*_{\bC^*}(\bP^r) = H^*(\bP^r_{\bC^*})$, we recall the
general formula for cohomology of a projective bundle.
Let $E\to X$ be a rank $(r+1)$ complex vector bundle over a 
topological space $X$, and let $\pi: \bP(E)\to  X$ be 
the projectivization of $E$, which is an $\bP^r$-bundle over $X$.
The cohomology $H^*(\bP(E))$ of the
total space $\bP(E)$ is an $H^*(X)$-algebra
generated by $H$ with a single relation
$$
H^{r+1} + c_1(E)H^r +\cdots + c_{r+1}(E) =0,
$$
where $c_i(E)$ is the $i$-th Chern class of $E$, and
$H$ is of degree 2.

In our case $E = \oplus_{i=0}^r \cO_{\bCP^\infty}(a_i)$, so 
the total Chern class of $E$ is given by
$$
c(E) = \prod_{i=0}^r (1+ a_i u).
$$
We have
$$
H^*_{\bC^*}(\bP^r)= 
H^*(\bP^r_{\bC^*}) \cong \bQ[u,H]/\langle 
\prod_{i=0}^r (H+ a_i u) \rangle,
$$
where $\bQ[u,H]$ is the ring of polynomials
in two variables $u$, $H$ with coefficients in $\bQ$, 
and $\langle \prod_{i=0}^r (H+ a_i u)\rangle$ is 
the principal ideal generated by $\prod_{i=0}^r (H+a_i u)$.
\end{exam}

\subsection{Equivariant vector bundle}

A continuous map $f:X\to Y$ between $G$-spaces
is called {\em $G$-equivariant} if 
$f(g\cdot x) = g\cdot f(x)$ for all $g\in G$ and $x\in X$.

Let $p: V\to X$ be a complex vector bundle over a
$G$-space $X$. We say $p:V\to X$ is a $G$-equivariant
complex vector bundle over $X$ if the following properties
hold.
\begin{itemize}
\item $V$ is a $G$-space.
\item $p$ is $G$-equivariant. 
\item For every $g\in G$, define $\tilde{\phi}_g: V\to V$ by $v\mapsto g\cdot v$,
and $\phi_g: X\to X$ by $x\mapsto g\cdot x$. Then 
$\tilde{\phi}_g$ is a vector bundle map covering $\phi_g$:
$$
\begin{CD}
V  @>{\tilde{\phi}_g}>>  V\\
@V{p}VV    @V{p}VV\\
X  @>{\phi_g}>>   X
\end{CD}
$$
\end{itemize}

\begin{exam}
When $X$ is a point, a complex vector bundle
over $X$ is a complex vector space $V$, and 
a $G$-equivariant vector bundle over $X$ is
a representation $\rho: G\to GL(V)$.
\end{exam}

\subsection{Equivariant Chern classes}
Let $\pi:V\to X$ be a $G$-equivariant vector bundle
over a $G$-space $X$. Then $V_G$ is a complex
vector bundle over $X_G$. The
$k$-th $G$-equivariant Chern class of $V$ is defined
to be the $k$-th Chern class of $V_G$:
$$
(c_k)_G(V) := c_k(V_G)\in H^{2k}(X_G)= H^{2k}_G(X).
$$
The $G$-equivariant Chern characters are defined
similarly:
$$
(\ch_k)_G(V) := \ch_k(V_G) \in H^{2k}(X_G)= H^{2k}_G(X).
$$
The $G$-equivariant Euler class of $V$ is defined to be
the Euler class (i.e. top Chern class) of $V_G$:
$$
e_G(V) := e(V_G) = c_r(V_G) \in H^{2r}(X_G)= H^{2r}_G(X)
$$
where $r=\rank_\bC V$.

\begin{exam}\label{Ca}
For any $a\in \bZ$, let
$\bC_a$ be the 1-dimensional representation of $\bC^*$ 
with character $t\mapsto t^a$. Then $\bC_a$ can be viewed
as a $\bC^*$-equivariant vector bundle over a point. We have
$$
(\bC_a)_{\bC^*} = \{ (u,v)\in (\bC^\infty-\{0\})\times \bC\}
/ t\cdot(u,v) \sim (t^{-1} u, t^a v)
\cong \cO_{\bCP^\infty}(-a)
$$
$$
(c_1)_{\bC^*}(\bC_a) = -au\in H^2_{\bC^*} (\pt)=\bZ u.
$$

\end{exam}

\subsection{Atiyah-Bott localization formula}

Suppose that $T=(\bC^*)^n$ acts on a complex manifold $X$, and
that the fixed points set $X^T$ is a disjoint union of compact complex
submanifolds $Z_1,\ldots, Z_N$ of $X$. Then 
the normal bundle $N_j$ of $Z_j$ in $X$ is a $T$-equivariant
complex vector bundle over $Z_j$. The equivariant Euler class 
$$
e_T(N_j) \in H_T^*(Z_j)\cong H^*(Z_j)\otimes_\bQ  H^*(BT) = H^*(Z_j)\otimes \bQ[u_1,\ldots, u_n]
$$
is not a zero divisor. Let $R=\bQ(u_1,\ldots,u_n)$. Then
$e_T(N_j)$ is invertible in $H^*(Z_j)\otimes_\bQ \bQ(u_1,\ldots, u_n)$.

Recall that $\pi: X_T \to BT$ is a fibration with fiber $X$. 
When $X$ is compact, the fiber $\pi^{-1}(b)$ over a point
$b\in BT$ represents a homology class $f\in H_{2n}(X_T)$, 
where $n=\dim_\bC X$; this class is independent
of choice of $b\in BT$. There is an additive map 
$$
\int_{X_T} :H^*_T(X) \to H^*_T(\pt),
$$
known as ``integration along the fiber'' or
``push-forward to a point''; it sends $H^q_T(X)$ to
$H^{q-2n}_T(\pt)$. (Intuitively, 
$\int_{X_T}$ is given by  contraction with $f$.)
Similarly, we have maps 
$$
\int_{(Z_j)_T}:H^*_T(Z_j)\to H^*_T(\pt)
$$ 
which send $H^q_T(Z_j)$ to $H^{q-2\dim_\bC Z_j}(\pt)$.
Here we will not give the precise definition of $\int_{X_T}$, but
the maps $\int_{(Z_j)_T}$ can be described very explicitly, as follows.
Any $\alpha \in H^*_T(Z_j)\cong H^*(Z_j)\otimes_\bQ \bQ[u_1,\ldots,u_n]$ is of the form
$$
\alpha =\sum_{i=1}^m \alpha_i p_i 
$$
where $\alpha_i\in H^*(X)$ and $p_i\in \bQ[u_1,\ldots,u_n]$.
Then
$$
\int_{(Z_j)_T} \alpha =\sum_{i=1}^m
\langle [Z_j], \alpha_i\rangle p_i
$$
where $[Z_j]\in H_{2\dim_\bC Z_j}(Z_j;\bQ)$ is the fundamental class, and
$\langle\ , \ \rangle$ is  the pairing between the homology 
$H_*(Z_j)$ and the cohomology $H^*(Z_j)$. We extend
$\int_{(Z_j)_T}$ to 
$$
\int_{(Z_j)_T}: H^*_T(Z_j)\otimes_{\bQ}\bQ(u_1,\ldots,u_n)\to \bQ(u_1,\ldots,u_n)
$$
by taking $p_i\in \bQ(u_1,\ldots,u_n)$.

\begin{theo}[Atiyah-Bott localization formula \cite{AB} ]
$$
\int_{X_T} \alpha =\sum_j \int_{(Z_j)_T} \frac{i_j^*\alpha}{e_T (N_j)}
$$
where $i_j:Z_j\hookrightarrow X$ is the inclusion.
\end{theo}

Let $X$ be a compact complex manifold with a holomorphic $T$-action. 
The constant map $X\to \pt$ also induces an additive map between
equivariant $K$-theories:
$$
\pi_!: K_T(X)\to K_T(\pt),\quad
\cE\mapsto \sum_i (-1)^i H^i(X,\cE)
$$ 
where $\cE$ is a $T$-equivariant holomorphic vector bundle
over $X$, and $H^i(X,\cE)$ are the sheaf cohomology groups, which
are representations of $T$.

A representation of $T$ is determined by its $T$-equivariant Chern 
character $\ch_T$. We can compute $\ch_T (\pi_!\cE)$ by Grothendieck-Riemann-Roch
(GRR) theorem and the Atiyah-Bott localization formula. Applying
GRR to  the fibration $\pi: X_T \to BT$, we have
$$
\ch_T (\pi_! \cE)  = \int_{X_T} \ch_T(\cE)\td_T(TX)
$$
where $\td_T(TX)$ is the $T$-equivariant Todd  class of
the tangent bundle $TX$ of $X$. 
By localization, 
$$
 \int_{X_T} \ch_T(\cE)\td_T(TX)
=\sum_{j=1}^N \int_{(Z_j)_T}\frac{i_j^*\left( \ch_T(\cE) \td_T(T_X)\right)}{e_T(N_j)}. 
$$

We now specialize to the case where $Z_j$ are isolated points.
We write $p_1,\ldots, p_N$ instead of $Z_1,\ldots, Z_N$.
Let $m=\dim_\bC X$, and let
$$
x_{j,1},\ldots, x_{j,m} \in 
H^2_T(\pt) = \bigoplus_{i=1}^n \bQ u_i 
$$
be the weights of the $T$-action on the tangent
space $T_{p_j}X$ of $X$ at $p_j$.   
Then
$$
i_j^* \td_T(T_X) =\prod_{k=1}^m \frac{x_{j,k}}{1-e^{-x_{j,k}}},\quad
e_T(N_j)= e_T(T_{p_j}X) = \prod_{k=1}^m x_{j,k}.
$$
Let $r=\rank_\bC \cE$, and let
$$
y_{j,1},\ldots, y_{j,r} \in 
H^2_T(\pt) 
$$
be the weights of the $T$-action on the fiber $\cE_{p_j}$ 
of $\cE$ at $p_j$. Then
$$
i_j^* \ch_T(\cE) = \sum_{l=1}^r e^{y_{j,l} }.
$$
Therefore
\begin{equation}\label{eqn:GRR}
\ch_T(\pi_!\cE) =\sum_{j=1}^N \frac{\sum_{l=1}^r e^{y_{j,l}} }
{\prod_{k=1}^m (1-e^{-x_{j,k}}) }.
\end{equation}

\begin{exam}\label{Pone}
Let $T=\bC^*$ act on $\bP^1$ by 
$$
t\cdot [x,y] = [t x, y].
$$
The fixed points set consists of two
isolated points $ q^0 = [0,1]$ and $q^1 =[1,0]$.
We have
$$
e_T(T_{q^0}\bP^1)=u,\quad
e_T(T_{q^1}\bP^1) = -u.
$$

For an equivariant lifting of
the line bundle $\cO(k) \to \bP^1$, we have
$$
e_T(\cO(k)_{q^0}) = au,\quad
e_T(\cO(k)_{q^1}) = (a-k)u
$$
for some $a\in \bZ$. 

By \eqref{eqn:GRR},
\begin{eqnarray*}
&& \ch_T H^0(\bP^1,\cO(k)) - \ch_T H^1(\bP^1, \cO(k))\\
&=&\frac{e^{au}}{1-e^{-u}} +\frac{e^{(a-k)u}}{1-e^u} 
=\frac{e^{au}(1-e^{-(k+1)u})}{1-e^{-u}}\\
&=& \begin{cases}
\sum_{i=0}^k e^{(a-i)u}, & k\geq 0\\
-\sum_{i=1}^{-k-1} e^{(a+i)u}, & k<0
\end{cases}
\end{eqnarray*}
Indeed we have
\begin{eqnarray*}
\ch_T H^0(\bP^1,\cO(k)) &=&
\begin{cases}
\sum_{i=0}^k e^{(a-i)u}, & k\geq 0\\
0, & k<0
\end{cases} \\
\ch_T H^1(\bP^1,\cO(k)) &=&
\begin{cases}
0, & k\geq -1\\
\sum_{i=1}^{-k-1} e^{(a+i)u}, & k <-1
\end{cases} 
\end{eqnarray*}

Given $a\in \bZ$,  let $\bC_a$ denote the irreducible
representation of $T=\bC^*$ characterized by
$e_T(\bC_a)= au$ (c.f.
Example \ref{Ca}). Then
\begin{eqnarray*}
H^0(\bP^1,\cO(k)) &=&
\begin{cases}
\oplus_{i=0}^k \bC_{a-i}, & k\geq 0\\
0, & k<0,
\end{cases} \\
H^1(\bP^1,\cO(k)) &=&
\begin{cases}
0, & k\geq -1,\\
\oplus_{i=1}^{-k-1} \bC_{a+i}, & k <-1.
\end{cases} 
\end{eqnarray*}
\end{exam}

\begin{exam}\label{orbiPone}
This example will arise in the localization calculations
in Section \ref{sec:localize}.

Let $T=\bC^*$ act on $\bP^1$ and on $\cO(k)$ as in
the previous example. Let $f: C\cong \bP^1\to \bP^1$ be the
degree $d$ map given by $[u,v]\mapsto [u^d,v^d]$.
Let $p^0=[0,1], p^1=[1,0]\in C$. Then
$f(p^i)=q^i$ for $i=0,1$, and $f^*T\bP^1 \cong (TC)^{\otimes d}$.
We have
\begin{eqnarray*}
e_T(T_{p^0}C)&=& \frac{1}{d}e_T(T_{q^0}\bP^1)=\frac{u}{d},\\
e_T(T_{p^1}C)&=& \frac{1}{d}e_T(T_{q^1}\bP^1)=-\frac{u}{d}, \\
e_T((f^*\cO(k))|_{p^0}) &=& e_T(\cO(k)|_{q^0})= au,\\
e_T((f^*\cO(k))|_{p^1}) &=& e_T(\cO(k)|_{q^1})=(a-k)u.
\end{eqnarray*}

By \eqref{eqn:GRR},
\begin{eqnarray*}
&& \ch_T H^0(C,f^*\cO(k)) - \ch_T H^1(C, f^*\cO(k))\\
&=&\frac{e^{au}}{1-e^{-u/d}} +\frac{e^{(a-k)u}}{1-e^{u/d}} 
=\frac{e^{au}(1-e^{-(kd+1)u/d})}{1-e^{-u/d}}\\
&=& \begin{cases}
\sum_{i=0}^{kd} e^{au-\frac{iu}{d} }, & k\geq 0\\
-\sum_{i=1}^{-kd-1} e^{au+\frac{iu}{d}}, & k<0
\end{cases}
\end{eqnarray*}
Indeed we have
\begin{eqnarray*}
\ch_T H^0(C,f^*\cO(k)) &=&
\begin{cases}
\sum_{i=0}^{kd} e^{ au -\frac{iu}{d}}, & k\geq 0\\
0, & k<0
\end{cases} \\
\ch_T H^1(C,f^*\cO(k)) &=&
\begin{cases}
0, & k\geq 0\\
\sum_{i=1}^{-kd-1} e^{au+\frac{iu}{d}}, & k < 0
\end{cases} 
\end{eqnarray*}

Given $a\in \bZ$, let $\bC_{\frac{a}{d}}$ denote the orbibundle
on $\mathrm{pt}/\bZ_d$ characterized by $e_T(\bC_{\frac{a}{d}})=\frac{au}{d}$. Then
\begin{eqnarray*}
H^0(C,f^*\cO(k)) &=&
\begin{cases}
\oplus_{i=0}^{kd} \bC_{a-\frac{i}{d}}, & k\geq 0\\
0, & k<0,
\end{cases} \\
H^1(C,f^*\cO(k)) &=&
\begin{cases}
0, & k\geq 0,\\
\oplus_{i=1}^{-kd-1} \bC_{a+\frac{i}{d}}, & k < 0.
\end{cases} 
\end{eqnarray*}
\end{exam}

\section{Proof of the ELSV Formula by Virtual Localization}
\label{sec:proof}

In this section, we present the proof of the ELSV formula
by virtual localization on moduli spaces of relative
stable maps, following Graber-Vakil \cite{GV1}.
In Section \ref{sec:relative-moduli},
we define moduli spaces of relative stable maps
$\Mbar_{g,0}(\bP^1,\mu)$. In Section \ref{sec:branch},
we identify the Hurwitz number $H_{g,\mu}$ with
a top intersection on $\Mbar_{g,0}(\bP^1,\mu)$, so
that it can be viewed as a relative Gromov-Witten
invariant for the pair $(\bP^1,\infty)$. 
Section \ref{sec:localize} contains the localization
calculations (using the torus action introduced
in Section \ref{sec:T-action}) which yield
the ELSV formula. 

\subsection{Moduli spaces} \label{sec:relative-moduli}
We fix a pair $(g,\mu)$, where $g$ is a nonnegative integer
(which will be the genus) and  $\mu=(\mu_1\geq \cdots \geq\mu_h>0)$ 
is a partition. 
Let $\cM_{g,0}(\bP^1,\mu)$ be the moduli space 
of ramified covers
$$
f: (C,x_1,\ldots, x_h) \to (\bP^1,q^1)
$$
of degree $d$ from a nonsingular complex algebraic curve (Riemann surface) $C$ of genus $g$ to 
$\bP^1$ such that the ramification type over a distinguished point $q^1=\infty \in \bP^1$ is specified
by the partition $\mu$, i.e.,
$$
f^{-1}(q^1)=\mu_1 x_1+\cdots +\mu_h x_h
$$ 
as Cartier divisors.  The moduli space $\cM_{g,0}(\bP^1,\mu)$ is not compact. To compactify it, we
consider the moduli space
$$
\Mbar_{g,0}(\bP^1,\mu)
$$
of relative stable maps to $(\bP^1,q^1)$. 
Moduli spaces of relative stable maps were first constructed
in the symplectic category by  A.Li-Ruan \cite{LR} and by Ionel-Parker \cite{IP1, IP2};
later J. Li constructed such moduli spaces in the algebraic category \cite{Li1, Li2}.
We need to use the algebraic version constructed by J. Li because virtual localization
on moduli spaces of relative stable maps has only been proved in the algebraic setting
\cite{GP, GV2}.
Moduli spaces of relative stable maps  are defined for a general pair $(X,D)$ where
$X$ is a smooth projective variety, and $D$ is a smooth divisor in $X$. Here
we content ourselves with the definition for the special case we need for
the ELSV formula: $X=\bP^1$ and $D$ is a point.

We use notation similar to that in \cite{LLZ1} and \cite{Liu}.
 Let $\bP^1(m)=\bP^1_1\cup \cdots \cup \bP^1_m$
be a chain of $m$ copies of $\bP^1$. For $l=1,\ldots,m-1$, let 
$q^1_l$ be the node at which $\bP^1_l$ and $\bP^1_{l+1}$ intersect. Let
$q_0^1 \in \bP^1_l$ and $q_m^1 \in \bP^1_m$ be smooth points.

A point of $\Mbar_{g,0}(\bP^1,\mu)$ is a morphism
$$
f:(C,x_1,\ldots,x_h)\to (\bP^1[m], q^1_m)
$$
where $C$ is complex algebraic curve of (arithmetic) genus $g$, 
with at most nodal singularities, and $\bP^1[m]$ is obtained by identifying
$q^1\in \bP^1$ with $q^1_0\in \bP^1(m)$. In particular, $\bP^1[0]=\bP^1$.  
We call the original $\bP^1=\bP^1_0$ the {\em root} component and $\bP^1_1,\ldots,\bP^1_m$
the {\em bubble} components. For $l=0,\ldots,m$, let
$C_l=f^{-1}(\bP^1_l)$, so that $C=C_0\cup C_1 \cup \cdots \cup C_m$, and
let $f_l: C_l\to \bP^1_l$ be the restriction of $f$. Then $f$ satisfies the following
properties:
\begin{enumerate}
\item ({\em degree}) $\deg f_l=d$, for $l=0,\ldots, m$.
\item ({\em ramification}) $f^{-1}(q^1_m)=\sum_{i=1}^h \mu_i x_i$ as Cartier divisors.
\item ({\em predeformability}) The preimage of each node of the target consists of
nodes, at each of which two branches have the same contact order;
distinct $C_i$ share no common irreducible components. This is the predeformable
condition: so that one can smooth both the target and the domain to obtain a morphism
to $\bP^1$.
\item ({\em stability}) The automorphism group of $f$ is finite.
\end{enumerate}

Two morphisms satisfying (1)--(3) are equivalent
if (a) they have the same target $\bP^1[m]$ for some nonnegative
integer $m$,  and  (b) they differ by an isomorphism of the
domain and an element of $\Aut(\bP^1(m),q^1_0, q^1_m)\cong (\bC^*)^m$. 
In particular, this defines the automorphism group in (4).
For fixed $g,\mu$, the stability condition (4) gives an upper bound of the number
$m$ of bubble components of the target.

By the results in  \cite{Li1}, 
$\Mbar_{g,0}(\bP^1,\mu)$ is a proper Deligne-Mumford stack
with a perfect obstruction theory of virtual dimension $r:= 2g-2+d+h$.
Roughly speaking, this means that it is a compact, Hausdorff, singular
orbifold,  with a ``virtual tangent bundle'' of rank $r$ (over $\bC$).  
It has a virtual fundamental class
$$
[\Mbar_{g,0}(\bP^1,\mu)]^{\vir} \in H_{2r}
(\Mbar_{g,0}(\bP^1,\mu);\bQ).
$$

\subsection{Branch morphism}\label{sec:branch}
Given a point
$$
[f:(C,x_1,\ldots, x_h)\to \bP^1[m] ] \in \Mbar_{g,0}(\bP^1,\mu),
$$
let $\tilde{f}= \pi_m\circ f:(C_1,x_1,\ldots,x_h)\to \bP^1$, where
$\pi_m:\bP^1[m]\to \bP^1$ is the map that contracts the bubble components.
When $C$ is nonsingular (in this case we must have $m=0$ and $\tilde{f}=f$), 
let $\Br(\tilde{f})$ denote the branch divisor of $\tilde{f}$, namely, 
$$
\Br(\tilde{f})=\sum_{q\in\bP^1}   (d-\# \tilde{f}^{-1}(q))q
$$
where $\# \tilde{f}^{-1}(q)\leq d$, and $\# \tilde{f}^{-1}(q) <d$
if and only if $q$ is a critical value of $\tilde{f}$. 
Then 
$$
\Br(\tilde{f}) = b_1+\cdots + b_r + (d-h)q^1
$$
where $b_1,\ldots, b_r\in \bP^1-\{q^1\}$ (not necessarily
distinct). In general $\Br(\tilde{f})= \Br'(\tilde{f}) + (d-h)q^1$, where
$\Br'(\tilde{f})$ is an effective divisor on $\bP^1$ of degree $r$
and can be viewed as a point in $\Sym^r \bP^1$, the $r$-th symmetric
produce of $\bP^1$; the map
\begin{eqnarray*}
\Br:\Mbar_{g,0}(\bP^1,\mu) &\to& \Sym^r \bP^1\cong \bP^r\\
{ [f:(C,x_1,\ldots, x_h)\to \bP^1[m] ] } &\mapsto& 
\Br'(\tilde{f})
\end{eqnarray*}
is a morphism, known as the {\em branch morphism} (c.f.
Fantechi-Pandharipande \cite{FanP}).
The Hurwitz number $H_{g,\mu}$ is essentially the degree
of the branch morphism $\Br:\Mbar_{g,0}(\bP^1,\mu)\to \bP^r$. More precisely,
let $H \in H^2(\bP^r;\bZ)$ be the hyperplane class, so
that $H^*(\bP^r;\bZ)=\bZ[H]/\langle H^{r+1} \rangle$. Then
$H^r\in H^{2r}(\bP^r;\bZ)$ is the Poincar\'{e} dual of the
point class $[\pt]\in H_0(\bP^r;\bZ)$, and 
$H^{2i}(\bP^r;\bQ)=\bQ H^i$, $i=1,\ldots,r$. We have
\begin{equation}\label{eqn:degree}
H_{g,\mu}=\frac{1}{\#\Aut(\mu)} \deg \Br =\frac{1}{\#\Aut(\mu)}
\int_{[\Mbar_{g,0}(\bP^1,\mu)]^\vir} \Br^* H^r
\end{equation}
where
$$
[\Mbar_{g,0}(\bP^1;\mu)]^\vir \in H_{2r}\left(\Mbar_{g,0}(\bP^1\mu);\bQ\right),\quad
\Br^* H^r \in H^{2r}\left(\Mbar_{g,0}(\bP^1,\mu);\bQ\right).
$$

We will use virtual localization to compute the right hand side
of \eqref{eqn:degree}  and obtain the right hand side
of the ELSV formula \eqref{eqn:ELSV}.

\subsection{Torus action}\label{sec:T-action}
Let $\bC^*$ act on $\bP^1$ by
$t\cdot [x,y]= [t x, y]$ where $t\in \bC^*$ and $[x,y]\in \bP^1$. Then
$\bC^*$ acts on $\Mbar_{g,0}(\bP^1,\mu)$ and on $\Sym^r\bP^1\cong \bP^r$, and the branch morphism
$\Br:\Mbar_{g,0}(\bP^1,\mu)\to \bP^r$ is $\bC^*$-equivariant. 
The isomorphism $\bP^r\stackrel{\cong}{\longrightarrow} \Sym^r(\bP^1)$ is given by
$$
[a_0, a_1,\ldots, a_r] \mapsto \mathrm{div}(\sum_{i=0}^r a_i x^i y^{r-i})
$$
where $[x,y]$ are homogeneous coordinates of $\bP^1$, and
$\mathrm{div}(p(x,y))$ is the divisor defined by the equation
$p(x,y)=0$.
The $\bC^*$-action on $\bP^r$ is given by
$t\cdot[a_0,\ldots,a_r]=[ a_0, t^{-1}a_1,\ldots, t^{-r} a_r]$
for $t\in \bC^*$, $[a_0,\ldots,a_r]\in \bP^r$. By Example \ref{Pr},
the $\bC^*$-equivariant cohomology of $\bP^r$ is given by
$$
H^*(\bP^r;\bQ) =\bQ[H,u]/\langle H(H-u)\cdots (H-ru)\rangle.
$$
The $\bC^*$-fixed points on $\bP^r$ are
$$
p_i= \mathrm{div}(x^i y^{r-i})= iq^0 + (r-i) q^1,\quad i=0,\ldots,r. 
$$

\subsection{Localization}\label{sec:localize}
We lift $H^r \in H^{2r}(\bP^r;\bQ)$, the Poincar\'{e} dual
of the point class $[\mathrm{pt}]\in H_0(\bP^r;\bQ)$, to
$\prod_{i=0}^{r-1}(H-iu) \in H^{2r}_{\bC^*}(\bP^r;\bQ)$,
the $\bC^*$-equivariant Poincar\'{e} dual of  the fixed
point $p_r \in \bP^r$. Then
\begin{eqnarray*}
H_{g,\mu} &=&\frac{1}{|\Aut(\mu)|} \int_{[\Mbar_{g,0}(\bP^1,\mu)]^{\vir}_{\bC^*} }
\Br^*\prod_{i=0}^{r-1}(H-iu)\\
&=& \frac{1}{|\Aut(\mu)|}\sum_F  \int_{[F]^{\vir}} 
\frac{\left(\Br^*\prod_{i=0}^{r-1}(H-iu)\right)\Bigr|_F}{e_{\bC^*}(N_F^\vir)}
\end{eqnarray*}
where the sum is over all connected components of
the $\bC^*$ fixed point set $\Mbar_{g,0}(\bP^1,\mu)^{\bC^*}$, and
$e_{\bC^*}(N_F^\vir)$ is the $\bC^*$-equivariant
Euler class of the virtual normal bundle $N_F^\vir$ of
the fixed locus $F$.  If $\Mbar_{g,0}(\bP^1,\mu)$ were
a compact complex manifold, and each $F$ were a compact
complex submanifold of dimension $d_F$  then $[F]^\vir$ would the usual
fundamental class $[F]\in H_{2d_F}(F;\bQ)$, 
$N^\vir_F$ would be the usual normal bundle $N_F$ in $\Mbar_{g,0}(\bP^1,\mu)$,
and the second equality would be the Atiyah-Bott localization formula.
Here we need to use the virtual localization formula proved
by Graber-Pandharipande \cite{GP}. By \cite{GV2}, we may
apply virtual localization  to moduli space of
relative maps $\Mbar_{g,0}(\bP^1,\mu)$.

For each $F$, we have $\Br(F) = p_i$ for some $i\in \{0,1,\ldots,r\}$.
Note that $H|_{p_i}= iu$. Let 
$$
F_i= \Mbar_{g,0}(\bP^1,\mu)^{\bC^*}\cap \Br^{-1}(p_i), \quad i=0,\ldots, r.
$$
Then $\Mbar_{g,0}(\bP^1,\mu)^{\bC^*} =F_0\cup \cdots \cup F_r$, and
$$
\Br^*(\prod_{j=0}^{r-1}(H-ju))\Bigr|_{F_i} = \prod_{j=0}^{r-1}(H-ju)\Bigr|_{p_i}
=\begin{cases} 0, & 0\leq i\leq r-1,\\ r! u^r, & i=r.\end{cases}
$$
Therefore,
$$
H_{g,\mu}=\frac{1}{|\Aut(\mu)|} \int_{[F_r]^\vir}\frac{r! u^r}{e_{\bC^*}(N^\vir_{F_r})}. 
$$

We need to identify the substack $F_r$, and its virtual tangent bundle
$T^\vir_{F_r}$ and virtual normal bundle $N^\vir_{F_r}$. 

Let $f:(C,x_1,\ldots,x_h)\to(\bP^1[m],q^1_m)$ be a relative stable map
which represents a point in $F_r$. One can show that $f$ must be of the
following form.
\begin{enumerate}
\item $m=0$, so $C=C_0$.
\item $D_0: =f^{-1}(q^0)$ is a curve of arithmetic genus $g$.
\item $f^{-1}(q^1)=\{ x_1,\ldots, x_h\}$.
\item $f^{-1}(\bP^1-\{q^0,q^1\})$ is a disjoint union of of twice punctured spheres
$L_1,\ldots,L_h$, and $f|_{L_i}: L_i\to \bP^1-\{q^0, q^1\}$ is an honest
covering map of degree $\mu_i$.
\item For $i=1,\ldots, h$, let $D_i$ be the closure of $L_i$ in $C$.
Then $D_i =L_i\cup \{ y_i, x_i\} \cong \bP^1$, where $x_i$ is the $i$-th marked
point, and $y_i$ is a node at which $D_i$ intersects $D_0$. $\hat{f}_i:= f|_{D_i}:D_i\to \bP^1$
is fully ramified over $q^0$ and $q^1$.
\item $(D_0,y_1,\ldots,y_h)$ is an $h$ pointed, genus $g$ stable curve, so
it represents a point in $\Mbar_{g,h}$. 
\end{enumerate} 
We have 
$$
\Aut(f)=\prod_{i=1}^h \Aut\left(\hat{f}_i:(D_i, y_i,x_i)\to (\bP^1,q^0,q^1)\right) =\prod_{i=1}^h\bZ_{\mu_i}
$$
where $\bZ_{\mu_i}$ is the cyclic group of order $\mu_i$.
Note that $\hat{f}_i:(D_i,y_i,x_i) \to (\bP^1,q^0,q^1)$ are the same
over $F_r$ while $(D_0,y_1,\ldots,y_h)$ can be any point in $\Mbar_{g,h}$. We have
a surjective morphism
$$
\iota: \Mbar_{g,h} \to  F_r = \Mbar_{g,h}/\prod_{i=1}^h \bZ_{\mu_i} \cong\Mbar_{g,h}\times
\prod_{i=1}^h B\bZ_{\mu_i}.
$$

To simplify the notation, we write $(C,\vx)$ instead of $(C,x_1,\ldots,x_h)$,
and write $(C,\vx,f)$ instead of $f:(C,x_1,\ldots,x_n)\to \bP^1$.
Let $D_0, D_1,\ldots,D_h$ and $\hat{f}_i :=f|_{D_i}$ be defined as above,
and write $(D_0,\vy)$ instead of $(D_0,y_1,\ldots,y_h)$.

The tangent space $T^1_{(C,\vx,f)}$ and the obstruction space $T^2_{(C,\vx,f)}$ at
a moduli point $[(C,\vx,f)]\in F_r$ fit in
the following long exact sequence
of $\bC^*$-representations:
\begin{equation}\label{eqn:les}
\begin{aligned}
0& \to  \Aut(C,\vx) \to \Def(f) \to T^1_{(C,\vx, f)}\\
 & \to  \Def(C,\vx) \to \Obs(f) \to T^2_{(C,\vx, f)}\to 0.
\end{aligned}
\end{equation}
\begin{itemize}
\item $\Aut(C,\vx)=\Ext^0(\Omega_C(\sum_{i=1}^h x_i),\cO_C)$
is the space of infinitesimal automorphism of the domain $(C,\vx)$.
We have
$$
\Aut(C,\vx)=\Aut(D_0, \vy)\oplus
\bigoplus_{i=1}^h \Aut(D_i,y_i,x_i)
$$
where $\Aut(D_0,\vy)=0$
since $(D_0,\vy)$ is stable, and
$\Aut(D_i, y_i, x_i) =H^0(\bP^1,T_{\bP^1}(-0-\infty))\cong \bC_0$
(trivial 1-dimensional representation of $\bC^*$).

\item  $\Def(C,\vx)=\Ext^1(\Omega_C(-\sum_{i=1}^h x_i),\cO_C)$
is the space of infinitesimal deformation of the domain $(C,\vx)$.
We have a short exact sequence of $\bC^*$-representations:
$$
0\to \Def(D_0,\vy)\to \Def(C,\vx)\to 
\bigoplus_{i=1}^h T_{y_i}D_0\otimes T_{y_i}D_i\to 0
$$
where $\Def(D_0,\vy)= T_{(D_0,\vy)}\Mbar_{g,h}$.

\item $\Def(f)=H^0\left(C,f^*(T\bP^1(-q^1))\right)$ 
is the space of infinitesimal deformation of the map $f$, and

\item $\Obs(f)= H^1\left(C,f^*(T\bP^1(-q^1))\right)$ is the space
of obstruction to deforming $f$. 
\end{itemize}

For $i=1,2$, let $T^{i,f}$ and $T^{i,m}$ be the fixed
and moving parts of $T^i\bigr|_{F_r}$. Then
$$
T^1 = T^{1,f} + T^{2,m},\quad
T^2 = T^{2,f} + T^{2,m}.
$$
The virtual tangent bundle of $F_r$ is $T_{F_r}^\vir = T^{1,f}-T^{2,f}$
and the virtual normal bundle of $F$ is $N_{F_r}^\vir = T^{1,m}-T^{2,m}$.
Let 
$$
B_1 =\Aut(C,\vx),\quad B_2=\Def(f),\quad B_4 = \Def(C,\vx),
\quad B_5 =\Obs(f),
$$
and let $B_i^f$ and $B_i^m$ be the fixed and moving parts
of $B_i$. Then 
$$
T_{F_r}^\vir = T^{1,f}-T^{2,f},
$$
$$
\frac{1}{e_{\bC^*}(N^\vir_{F_r})}=\frac{e_{\bC^*}(B_5^m)e_{\bC^*}(B_1^m)}{e_{\bC^*}(B_2^m)e_{\bC^*}(B_4^m)}
$$
We have
\begin{eqnarray*}
&& B_1^f =\bC_0^{\oplus h},\quad B_1^m=0,\\
&& B_4^f = T_{(D_0,\vy)} \Mbar_{g,h},\quad B_4^m =\oplus_{i=1}^h T_{y_i}D_0 \otimes T_{y_i}D_i
=\oplus_{i=1}^h (\bL_i^\vee)_{(C_0,\vy)} \otimes \bC_{\frac{1}{\mu_i}}.
\end{eqnarray*}
Let $\iota^*: H^*(F_r)\otimes \bQ(u) \to H^*(\Mbar_{g,h})\otimes_\bQ \bQ(u)$
be induced by the surjective morphism $i:\Mbar_{g,h}\to F_r$. Then
$$
\iota^*\left(\frac{e_{\bC^*}(B_1^m)}{e_{\bC^*}(B_4^m)}\right)
=\frac{1}{\prod_{i=1}^h (\frac{u}{\mu_i}-\psi_i)} 
=\frac{\mu_1\cdots \mu_h}{\prod_{i=1}^h(u-\mu_i\psi_i)}.
$$

For $k=0,1$ and $i=0,1,\ldots,h$, let
$$
H^k(D_i) = H^k\left(D_i, \hat{f}_i^*(T\bP^1(-q^1))\right).
$$
Then we have a long exact sequence of $\bC^*$-representations:
\begin{eqnarray*}
0&\to & B_2\to H^0(D_0)\oplus \bigoplus_{i=1}^h H^0(D_i) \to (T_{q^0}\bP^1)^{\oplus h}\\
&\to & B_5 \to H^1(D_0)\oplus \bigoplus_{i=1}^h H^1(D_i)\to 0.
\end{eqnarray*}
We have
$$
H^k(D_0)  \cong H^k(D_0)\otimes T_{q^0}\bP^1
=\begin{cases} \bC_{1}, & k=0;\\
(\bE^\vee)_{(D_0,\vy)} \otimes \bC_{1}, & k=1.
\end{cases} 
$$
Note that
$$
T\bP^1(-q^1)\cong \cO_{\bP^1}(1),
\quad e_{\bC^*}\left ((T\bP^1(-q^1))_{q^0}\right) = u.
$$
By Example \ref{orbiPone}, for $i=1,\ldots,h$, we have
$$
H^k(D_i)\cong H^k(D_i,\hat{f}_i^*\cO(1))
= \begin{cases} 
\oplus_{a=0}^{\mu_i} \bC_{\frac{a}{\mu_i}}, & k=0;\\
0, & k=1.
\end{cases}
$$
Therefore, 
$$
B_2^f=\bC_0^{\oplus h},\quad B_5^f =0,
\quad B_5^m -B_2^m = \bE^\vee_{(D_0,\vy)}\otimes\bC_1 \oplus \bC_1^{\oplus (h-1)}
-\bigoplus_{i=1}^h\bigoplus_{a=1}^{\mu_i}\bC_{\frac{a}{\mu_i}}. 
$$
So 
$$
\iota^*\left(\frac{e_{\bC^*}(B_5^m)}{e_{\bC^*}(B_2^m)}\right)=
\prod_{i=1}^h\frac{\mu_i^{\mu_i}}{\mu_i!} u^{h-d-1}\cdot \Lambda^\vee_g(u),
$$
where 
$$
\Lambda_g^\vee(u)=\sum_{i=0}^g (-1)^i\lambda_i u^{g-i}.
$$
We conclude that
$$
T^{1,f}_{(C,\vx,f)}= T_{(D_0,\vy)}\Mbar_{g,h},\quad T^{2,f}=0, 
$$
so
$$
[F_r]^\vir = \frac{1}{\mu_1\cdots\mu_h} \iota_*[\Mbar_{g,h}]
$$
where $\iota_*:H_{2r}(\Mbar_{g,h};\bQ) \to H_{2r}(F_r;\bQ)$. We also have
$$
\frac{1}{\iota^*e_{\bC^*}(N^\vir_{F_r})}=\prod_{i=1}\frac{\mu_i^{\mu_i}}{\mu_i!}
\cdot \frac{\mu_1\cdots \mu_h \Lambda^\vee_g(u)u^{h-d-1}}{\prod_{i=1}^h(u-\mu_i \psi_i)}.
$$

\begin{eqnarray*}
H_{g,\mu}&=&\frac{1}{\#\Aut(\mu)}\int_{[F_r]^\vir} \frac{r! u^r}{e_{\bC^*} (N^\vir_{F_r})}\\
&=&\frac{r!}{\mu_1\cdots \mu_h \cdot \#\Aut(\mu)}
\int_{\Mbar_{g,h}}\frac{u^r}{\iota^* e_{\bC^*}(N_{F_r}^\vir)}\\
&=& \frac{r!}{\# \Aut(\mu)}\prod_{i=1}^h\frac{\mu_i^{\mu_i}}{\mu_i !}
\int_{\Mbar_{g,h}}\frac{\Lambda_g^\vee(u)  u^{2g-3+2h}}
{\prod_{i=1}^h(u-\mu_i\psi_i)}\\
&\stackrel{u=1}{=}& \frac{r!}{\# \Aut(\mu)}\prod_{i=1}^h\frac{\mu_i^{\mu_i} }{\mu_i !}
\int_{\Mbar_{g,h}}\frac{\Lambda_g^\vee(1)}
{\prod_{i=1}^h(1-\mu_i\psi_i)}.
\end{eqnarray*}
In the last equality,  we may let $u=a$, where $a$ is any nonzero constant rational
number; the answer is independent of $a$.
This completes the proof of the ELSV formula \eqref{eqn:ELSV}.


\begin{thebibliography}{AA}

\bibitem{AB}  M. F. Atiyah and R. Bott, 
``The moment map and equivariant cohomology,''
Topology {\bf 23} (1984), no. 1, 1--28. 

\bibitem{BF} K. Behrend and B. Fantechi,
``The intrinsic normal cone,''
Invent. Math. {\bf 128} (1997), 45--88.

\bibitem{CLL} L. Chen, Y. Li, and K. Liu,
``Localization, Hurwitz numbers and the Witten conjecture,'' 
Asian J. Math. {\bf 12} (2008), no. 4, 511--518. 

\bibitem{CK} D.A. Cox and S. Katz, 
{\it Mirror symmetry and algebraic geometry},
Mathematical Surveys and Monographs {\bf  68}, American Mathematical Society, Providence, RI, 1999.

\bibitem{DM} P. Deligne and D. Mumford,
``The irreducibility of the space of curves of given genus,''
Inst. Hautes \'{E}tudes Sci. Publ. Math. No. 36 (1969), 75--109. 
 
\bibitem{ELSV1} T. Ekedahl, S. Lando, M. Shapiro, and A. Vainshtein, 
``On Hurwitz numbers and Hodge integrals,''
C. R. Acad. Sci. Paris S\'{e}r. I Math. {\bf 328} (1999), no. 12, 1175--1180. 

\bibitem{ELSV2} T. Ekedahl, S. Lando, M. Shapiro, and A. Vainshtein,
``Hurwitz numbers and intersections on moduli spaces of curves,''
Invent. Math. {\bf 146} (2001), no. 2, 297--327.

\bibitem{Fa} C. Faber, 
``Algorithms for computing intersection numbers on moduli spaces of curves, 
with an application to the class of the locus of Jacobians,''
New trends in algebraic geometry (Warwick, 1996), 93--109, 
London Math. Soc. Lecture Note Ser., {\bf 264}, Cambridge Univ. Press, Cambridge, 1999. 

\bibitem{FP} C. Faber and R. Pandharipande,
``Hodge integrals and Gromov-Witten theory,''
Invent. Math. {\bf 139} (2000), no. 1, 173--199.

\bibitem{FanP} B. Fantechi and R. Pandharipande,
``Stable maps and branch divisors,''
Compositio Math. {\bf 130} (2002), no. 3, 345--364. 

\bibitem{Fu} W. Fulton, 
``Equivariant Cohomology in Algebraic Geometry,'' 
Eilenberg lectures at Columbia University, Spring 2007, 
notes by Dave Anderson are available at
http://www.math.washington.edu/$\sim$dandersn/eilenberg/

 
\bibitem{GJV} I.P. Goulden, D.M. Jackson, A. Vainshtein,
``The number of ramified coverings of the sphere by the torus and surfaces 
of higher genera,'' Ann. Comb. {\bf 4} (2000), no. 1, 27--46.
 
\bibitem{GP}  T. Graber, R. Pandharipande, 
``Localization of virtual classes,''
Invent. Math. {\bf 135}  (1999),  no. 2, 487--518. 

\bibitem{GV1} T. Graber, R. Vakil, 
``Hodge integrals and Hurwitz numbers via virtual localization,''
Compositio Math. {\bf  135} (2003), no. 1, 25--36. 

\bibitem{GV2} T. Graber, R. Vakil,
``Relative virtual localization and vanishing of tautological classes on moduli spaces of curves,''
Duke Math. J. {\bf 130} (2005), no. 1, 1–37.

\bibitem{IP1} E.-N. Ionel, T. Parker, 
``Relative Gromov-Witten invariants,''
Ann. of Math. (2) {\bf 157}  (2003),  no. 1, 45--96.

\bibitem{IP2} E.-N. Ionel, T. Parker,
``The symplectic sum formula for Gromov-Witten invariants,''
 Ann. of Math. (2) {\bf 159} (2004), no. 3, 935--1025. 

\bibitem{Ka} M.E. Kazarian, 
``KP hierarchy for Hodge integrals,''
Adv. Math. {\bf 221} (2009), no. 1, 1--21. 

\bibitem{KaL} M.E. Kazarian, S. K. Lando,
``An algebro-geometric proof of Witten's conjecture,''
J. Amer. Math. Soc. {\bf 20} (2007), no. 4, 1079--1089.


\bibitem{KiL} Y-S. Kim, K. Liu,
``Virasoro constraints and Hurwitz numbers through asymptotic analysis,''
Pacific J. Math. {\bf 241} (2009), no.2, 275--284.

\bibitem{Kn2} F. Knudsen,
``The projectivity of the moduli space of stable curves. II. The stacks $M_{g,n}$,'' 
Math. Scand. {\bf 52} (1983), no. 2, 161--199.  

\bibitem{Kn3} F. Knudsen,
``The projectivity of the moduli space of stable curves. III. The line bundles on $M_{g,n}$, and a proof of 
the projectivity of $\overline M_{g,n}$ in characteristic $0$,''
Math. Scand. {\bf 52} (1983), no. 2, 200--212.  

\bibitem{KM} F. Knudsen and D. Mumford, 
``The projectivity of the moduli space of stable curves. I. Preliminaries on ``det'' and ``Div'','' 
Math. Scand. {\bf 39} (1976), no. 1, 19--55. 


\bibitem{Ko1} M. Kontsevich, 
``Intersection theory on the moduli space of curves and the matrix Airy function,''
Comm. Math. Phys. {\bf 147} (1992), no. 1, 1--23.

\bibitem{Ko2} M. Kontsevich,
``Enumeration of rational curves via torus actions,''
The moduli space of curves (Texel Island, 1994), 335--368, 
Progr. Math., {\bf 129}, Birkh\"{a}user Boston, Boston, MA, 1995. 


\bibitem{Li1} J. Li,
``Stable Morphisms to singular schemes and relative stable morphisms,''
J. Diff. Geom. {\bf 57} (2001), 509-578.

\bibitem{Li2} J. Li,
``Relative Gromov-Witten invariants and a degeneration
formula of Gromov-Witten invariants,''
J. Diff. Geom. {\bf 60} (2002), 199-293.

\bibitem{LT} J. Li and G. Tian,
`` Virtual moduli cycles and Gromov-Witten invariants of
algebraic varieties,''
J. Amer.  Math. Soc. {\bf 11} (1998), 119--174.

\bibitem{LR} A. Li, Y. Ruan,
``Symplectic surgery and Gromov-Witten invariants of Calabi-Yau 3-folds,''
Invent. Math. {\bf 145} (2001), no. 1, 151--218.

\bibitem{LZZ} A. Li, G. Zhao, Q. Zheng, 
``The number of ramified covering of a Riemann surface by Riemann surface,''
Comm. Math. Phys. {\bf  213}  (2000),  no. 3, 685--696. 

\bibitem{LLZ1} C.-C. M. Liu, K. Liu, J. Zhou,
`` A proof of a conjecture of Mari\~{n}o-Vafa on Hodge integrals,''
J. Differential Geom. {\bf 65}  (2003),  no. 2, 289--340. 

\bibitem{LLZ2} C.-C. M. Liu, K. Liu, J. Zhou, 
``A formula of two-partition Hodge integrals,''
J. Amer. Math. Soc. {\bf 20} (2007), no. 1, 149--184

\bibitem{Liu} C.-C. M. Liu, 
``Formulae of one-partition and two-partition Hodge integrals,"
{\em The interaction of finite-type and Gromov-Witten invariants (BIRS 2003)}, 105--128, 
Geom. Topol. Monogr., {\bf 8}, Geom. Topol. Publ., Coventry, 2006. 

\bibitem{Ma} I.G. MacDonald,
{\em Symmetric functions and Hall polynomials,}
2nd edition, Claredon Press, 1995.


\bibitem{MV} M. Mari\~{n}o,  C. Vafa,
``Framed knots at large N,'' in {\em Orbifolds in mathematics
and physics} (Madison, WI, 2001), 185-204,
Contemp. Math. {\bf 310}, Amer. Math. Soc., Providence, RI, 2002.

\bibitem{Mi} M. Mirzakhani,
``Weil-Petersson volumes and intersection theory on the moduli space of curves,'' 
J. Amer. Math. Soc. {\bf 20} (2007), no. 1, 1--23.

\bibitem{MuZ} M. Mulase, N. Zhang,
``Polynomial recursion formula for linear Hodge integrals,''
{\tt arXiv:0908.2267}.


\bibitem{Mu} D. Mumford,
``Towards an enumerative geometry of the moduli space of curves,''
Arithmetic and geometry, Vol. II, 271--328, Progr. Math., {\bf 36}, 
Birkh\"{a}user Boston, Boston, MA, 1983. 

\bibitem{OP1} A. Okounkov, R. Pandharipande,
``Gromov-Witten theory, Hurwitz numbers, and matrix models,''
Algebraic geometry---Seattle 2005. Part 1, 325--414, Proc. Sympos. 
Pure Math., {\bf 80}, Part 1, Amer. Math. Soc., Providence, RI, 2009. 

\bibitem{OP2} A. Okounkov, R. Pandharipande, 
``Hodge integrals and invariants of the unknot,'' 
Geom. Topol. {\bf 8} (2004), 675--699. 


\bibitem{Wi} E. Witten, 
``Two-dimensional gravity and intersection theory on moduli space,'' 
Surveys in differential geometry (Cambridge, MA, 1990), 243--310, 
Lehigh Univ., Bethlehem, PA, 1991. 

\bibitem{Zh1} J. Zhou,
``Hodge integrals, Hurwitz numbers, and Symmetric Groups,''
{\tt arXiv:math/0308024}.

\bibitem{Zh2} J. Zhou,
``A Conjecture on Hodge Integrals,''
{\tt arXiv:math/0310282}.



\end{thebibliography}
\end{document}